
\input amssym

\hsize=15,4truecm
\vsize=22.5truecm
\advance\voffset by 1truecm
\mathsurround=1pt

\def\chapter#1{\par\bigbreak \centerline{\bf #1}\medskip}

\def\section#1{\par\bigbreak {\bf #1}\nobreak\enspace}

\def\sqr#1#2{{\vcenter{\hrule height.#2pt
      \hbox{\vrule width.#2pt height#1pt \kern#1pt
         \vrule width.#2pt}
       \hrule height.#2pt}}}

\def\k{\kappa}
\def\o{\omega}

\def\l{\lambda}

\def\a{\alpha}
\def\b{\beta}

\def\g{\gamma}

\def\n{\eta}


\def\A{{\cal A}}
\def\B{{\cal B}}
\def\C{{\cal C}}
\def\D{{\cal D}}

\def\M{{\bf M}}


\def\th #1 #2. #3\par\par{\medbreak{\bf#1 #2.
\enspace}{\sl#3\par}\par\medbreak}
\def\rem #1 #2. #3\par{\medbreak{\bf #1 #2.
\enspace}{#3}\par\medbreak}
\def\proof{{\bf Proof}.\enspace}
\def\sqr#1#2{{\vcenter{\hrule height.#2pt
      \hbox{\vrule width.#2pt height#1pt \kern#1pt
         \vrule width.#2pt}
       \hrule height.#2pt}}}
\def\eop{\mathchoice\sqr34\sqr34\sqr{2.1}3\sqr{1.5}3}

                                                                     %
                                                                     %
\newdimen\refindent\newdimen\plusindent                              %
\newdimen\refskip\newdimen\tempindent                                %
\newdimen\extraindent                                                %
                                                                     %
                                                                     %
\def\ref#1 #2\par{\setbox0=\hbox{#1}\refindent=\wd0                  %
\plusindent=\refskip                                                 %
\extraindent=\refskip                                                %
\advance\extraindent by 30pt                                         %
\advance\plusindent by -\refindent\tempindent=\parindent 
\parindent=0pt\par\hangindent\extraindent 
{#1\hskip\plusindent #2}\parindent=\tempindent}                      %
\refskip=\parindent                                                  %
                                                                     %

\def\empty{\emptyset}

\def\raj{\restriction}

\def\da{\downarrow}

\def\nda{\mathrel{\lower0pt\hbox to 3pt{\kern3pt$\not$\hss}\downarrow}}
\def\nDa{\mathrel{\lower0pt\hbox to 3pt{\kern3pt$\not$\hss}\Downarrow}}
\def\nbot{\mathrel{\lower0pt\hbox to 4pt{\kern3pt$\not$\hss}\bot}}
\def\ekom{\mathrel{\lower3pt\hbox to 0pt{\kern3pt$\sim$\hss}\mapsto}}
\def\do{\triangleright}

\def\anR{\mathrel{\lower1pt\hbox to 2pt{\kern3pt$R$\hss}\not}}
\def\anoR{\mathrel{\lower1pt\hbox to 2pt{\kern3pt$\overline{R}$\hss}\not}}

\def\anRm{\mathrel{\lower1pt\hbox to 2pt{\kern3pt$R^{-}$\hss}\not}}

\def\ndda{\mathrel{\lower0pt\hbox to 1pt{\kern3pt$\not$\hss}\downdownarrows}}

\def\warrow{\mathrel{\lower0pt\hbox to 1pt{\kern3pt$^{w}$\hss}\rightarrow}}

\null
\vskip 1truecm
\centerline{\bf MAIN GAP FOR LOCALLY SATURATED ELEMENTARY}
\centerline{\bf SUBMODELS OF A HOMOGENEOUS STRUCTURE}
\vskip 1truecm
\centerline{Tapani Hyttinen and Saharon Shelah$^{*}$}
\vskip 1.5truecm

\chapter{Abstract}
\bigskip

We prove a main gap theorem for locally saturated submodels of
a homogeneous structure. We also study the number
of locally saturated models, which are not elementarily embeddable
to each other.

\vskip 1.5truecm

Through out this paper we assume that $\M$ is a homogeneous model
of similarity type (=language) $L$. We study elementary submodels
of $\M$. We use $\M$ as the monster model
is used in stability theory
and so we assume that the
cardinality of $\M$ is large enough
for all constructions we
do in this paper. In fact, as in [HS1],
we assume that $\vert\M\vert$
is strongly inaccessible. Alternatively we could
assume less about $\vert\M\vert$ and instead of studying all
elementary submodels of $\M$, we could study suitably small
ones. Also the assumption that $\M$ is homogeneous can be replaced
by the assumption that $\M$ is $\k$-homogeneous for $\k$ large
enough. Notice that by [Sh1], if $D$ is a stable finite
diagram, then $D$ has a monster model like $\M$.

We assume that the reader is familiar with [HS1] and use
its notions and results freely.

\th 0.1 Definition.

(i) Suppose $\M$ is stable.
We say that $\A$ is $s$-saturated if it is
$F^{\M}_{\l (\M )}$-saturated i.e.
for all $A\subseteq\A$
of power $<\l (\M )$ and $a$ there is $b\in\A$ such that
$t(b,A)=t(a,A)$.

(ii) We say that $\A$ is locally $F^{\M}_{\k}$-saturated if for all
finite $A\subseteq\A$ there is $F^{\M}_{\k}$-saturated
$\B$ such that $A\subseteq\B\subseteq\A$. If $\M$ is stable,
then we say that
$\A$ is $e$-saturated if it is locally $F^{\M}_{\l (\M )}$-saturated.

(iii) Suppose $\M$ is stable.
We say that $\A$ is strongly $F^{\M}_{\k}$-saturated if
for all $A\subseteq\A$ of power $<\k$ and $a$ there is
$b\in\A$ such that $b\ E^{m}_{min, A}\ a$. By
$a$-saturated we mean strongly $F^{\M}_{\k (\M )}$-saturated.

\vskip 1truecm

\noindent
$*$ Research supported by the United States-Israel
Binational Science Foundation. Publ. 676.

\vfill
\eject

\th 0.2 Lemma.

(i) Every $F^{\M}_{\k}$-saturated model is locally $F^{\M}_{\k}$-saturated
and so (assuming $\M$ is stable) every $s$-saturated model is $e$-saturated.

(ii) Suppose $\M$ is stable. Then every $e$-saturated model is
strongly $F^{\M}_{\o}$-saturated.

(iii) Suppose $\M$ is superstable. Then every $e$-saturated model
is $s$-saturated.

\proof (i) is trivial and (ii) follows immediately from
[HS1] Lemma 1.9 (iv). So we prove (iii): Assume $\A$ is
$e$-saturated. Notice that by
(ii), $\A$ is $a$-saturated. Let $A\subseteq\A$ be of power $<\l (\M )$
and $a$ arbitrary. We show that there is $b\in\A$ such that
$t(b,A)=t(a,A)$. Clearly we may assume that $a\cap\A=\empty$.

Choose finite $B\subseteq\A$ so that $a\da_{B}\A$.
Since $\A$ is $e$-saturated,
we can find $s$-saturated $\B$ such that $B\subseteq\B\subseteq\A$.
Since by [HS1] Lemma 1.9 (iii) $\B$ is strongly
$F^{\M}_{\l (\M )}$-saturated, we can find $a_{i}\in\B$, $i<\l (\M )$,
such that $a_{i}\ E^{m}_{min, B}\ a$ and
$a_{i}\da_{B}\cup_{j<i}a_{j}$. Let $I=\{ a_{i}\vert\ i<\l (\M )\}$.
For all $i<\k (\M )$, choose $b_{i}$ so that
$t(b_{i},\A )=t(a,\A )$ and $b_{i}\da_{\A}\cup_{j<i}b_{j}$.
Let $J=\{ b_{i}\vert\ i<\k (\M )\}$. By [HS1] Corollaries 3.5 (iv)
and 3.11, $I\cup J$ is indiscernible over $B$. So
$$Av(I,A)=Av(J,A)=t(a,A).$$
Since $\vert A\vert <\l (\M )$ is regular,
we can find $C\subseteq B\cup I$ of power $<\l (\M )$
such that for all $c\in A$, $t(c,B\cup I)$ does not
split strongly over $C$. Let $b\in I\ (\subseteq\B\subseteq\A )$ be such that
$b\cap C=\empty$. Then clearly $t(b,A)=Av(I,A)=t(a,A)$. $\eop$

We prove a main gap theorem for $e$-saturated submodels of $\M$.
At some extend, the proofs are similar to the related proofs
in the case of complete first-order theories. So some of
the proofs are sketchy.

\chapter{1. Regular types}

In (the end of) the next section, regular types are needed.
In this section we prove the basic properties and
the existence of regular types. In this section we assume that
$\M$ is stable.

\th 1.1 Definition.

(i) We say that a stationary pair $(p,A)$ is regular if
the following holds: if $C\supseteq dom(p)$, $a\models p$
and $a\nda_{A}C$, then $(p,A)$ is orthogonal to $t(a,C)$.

(ii) Assume $\A$ is $s$-saturated and $p\in S(\A )$.
We say that $p$ is regular, if there are $A\subseteq B\subseteq\A$
such that $p$ does not split strongly over $A$,
$(p\raj B,A)$ is a regular stationary pair and $\vert B\vert <\k (\M )$.

\th 1.2 Lemma. Assume $\A$ is $s$-saturated, regular $p\in S(\A )$
is not orthogonal to $t(a,\A )$ and $\B$ is $s$-primary over
$\A\cup a$. Then there is $b\in\B$ such that $t(b,\A )=p$.

\proof Assume not. Let $A\subseteq B\subseteq\A$ be as in Definition
1.1 (ii). For all $i<\k (\M )$ choose $\A_{i}$ as follows:

(i) $\A_{0}=\A$,

(ii) if $i$ is limit, then $\A_{i}\subseteq\B$ is $s$-primary over
$\cup_{j<i}\A_{j}$,

(iii) if $i=j+1$ and there is $b_{j}\in\B$ such that
$t(b_{j},B)=p\raj B$ and $a\nda_{\A_{j}}b_{j}$, then $\A_{i}\subseteq\B$ is
$s$-primary over $\A_{j}\cup b_{j}$, if such $b_{j}$ does not exist
then we let $\A_{i}=\A_{j}$.

\noindent
Clearly there is $i<\k (\M )$ such that $\A_{i}=\A_{i+1}$.
Let $i^{*}$ be the least such ordinal. Then

(*) $t(a,\A_{i^{*}})$ is orthogonal to $p$.

\noindent
Let $\A^{*}$ be $s$-primary over $\A_{i^{*}}\cup a$.

{\bf Claim.} Assume $b\models p$. Then
$p\vdash t(b,\A^{*})$.

\proof Since $p$ is not realized in $\B$,
for all $i<i^{*}$, $b_{i}\nda_{A}\A_{i}$
and so, since $p$ is regular,
for all $i<i^{*}$, $p$ is orthogonal to $t(b_{i},\A_{j})$.
By induction on $i\le i^{*}$ it is easy to see that
$p\vdash t(b,\A_{i^{*}})$. By (*) above,
$p\vdash t(b,\A^{*})$. $\eop$ Claim.

By Claim, $p$ is orthogonal to $t(a,\A )$, a contradiction. $\eop$

\th 1.3 Corollary. Assume $\A_{i}$, $i<3$, are $s$-saturated,
$p_{i}\in S(\A_{i})$ and $p_{1}$ is regular.
If $p_{0}$ is not orthogonal to $p_{1}$ and
$p_{1}$ is not orthogonal to $p_{2}$, then
$p_{0}$ is not orthogonal to $p_{2}$.

\proof Immediate by Lemma 1.2. $\eop$

\th 1.4 Lemma. Assume that $\A$ is $s$-saturated, $a\nda_{\A}b$ and
$t(b,\A )$ is regular. Then $a\do_{\A}b$.

\proof Let $\l =(\l (\M ))^{+}$.
Clearly we may assume that $\A$ is
$F^{\M}_{\l}$-saturated. For a contradiction, assume that there
is $c$ such that $c\da_{\A}a$ and $c\nda_{\A}b$.
Choose $A\subseteq B\subseteq\B\subseteq\A$ such that

(i) $(t(b,B),A)$ is a regular stationary pair and $b\da_{A}\A$,

(ii) $\vert B\vert <\k (\M )$ and $\vert\B\vert =\l (\M )$,

(iii) $\B$ is $s$-saturated and $a\cup b\cup c\da_{\B}\A$.

\noindent
Then $b\nda_{\B}a$, $b\nda_{\B}c$ ([HS1] Lemma 3.8 (iv))
and $a\da_{\B}c$. Let $\A^{*}$ be $F^{\M}_{\l}$-primary over
$\A\cup a$ and $\C\subseteq\A^{*}$ $s$-primary over $\B\cup a$.
Without loss of generality we may assume that
$b\cup c\da_{\C}\A$.

For all $i<\k (\M )$, choose $b_{i}\in\A^{*}$ such that
$t(b_{i},\C\cup\bigcup_{j<i}b_{j})=t(b,\C\cup\bigcup_{j<i}b_{j})$.
Let $I=\{ b_{i}\vert\ i<\k (\M )\}$. Then $I\cup\{ b\}$ is indiscernible
over $\C$. Since $b\nda_{\B}\C$, it is easy
to see that $I\cup\{ b\}$ is not $\B$-independent.
So we can choose finite $J\subseteq I$ such that

(*) $J\cup\{ b\}$ is not $\B$-independent.

\noindent
If $J$ is chosen
so that $\vert J\vert$ is minimal, then $J$ is $\B$-independent.

Let $\D$ be $s$-primary over $\B\cup c$. Then
by (iii), $J\da_{\B}\D$ and so
$J$ is $\D$-independent. Since $p$ is regular,
$J\da_{\D}b$ and so $J\da_{\B}b$. Clearly this contradicts (*) above.
$\eop$

Assume $\A$ is $s$-saturated and $a\not\in\A$. We write
$Dp(a,\A )>0$ if there is $s$-primary model $\B$ over
$\A\cup a$ and $b\not\in\B$ such that $t(b,\B )$ is orthogonal
to $\A$.

\th 1.5 Lemma. Assume that $\M$ is superstable without
$(\l (\M ))^{+}$-dop. Let $\A$ be $s$-saturated, $I$ be
$\A$-independent and $a\nda_{\A}I$. If $t(a,\A )$
is regular and $Dp(a,\A )>0$, then there is $b\in I$
such that $a\nda_{\A}b$. And so by Lemma 1.4,
$a\da_{\A}\cup (I-\{ b\})$.

\proof Assume not. Clearly we may assume that $\vert\A\vert =\l (\M )$.
Choose $a_{i}$, $\A_{i}$ and $\C_{i}$,
$i<\a^{*}$, so that

(i) $a\da_{\A}a_{i}$,

(ii) $\A_{i}$ is $s$-primary over $\A\cup a_{i}$,

(iii) $\{ a_{i}\vert\ i<\a^{*}\}$ is $\A$-independent,

(iv) $\C_{0}=\A_{0}$ and $\C_{i+1}$ is $s$-primary
over $\C_{i}\cup\A_{i+1}$,

(v) $a\nda_{\C_{i}}\A_{i+1}$,

(vi) $(a_{i})_{i<\a^{*}}$ is a maximal sequence satisfying
(i)-(v) above.

\noindent
Since $\M$ is superstable, $\a^{*}<\o$. Let $n$ be such that
$\a^{*}=n+1$.
Let $\l =(\l (\M ))^{+}$ and
$\B$ be $F^{\M}_{\l}$-saturated model such that
$\A\subseteq\B$ and $\B\da_{\A}\C_{n}$.
Let $\B_{i}$ be $F^{\M}_{\l}$-primary over $\B\cup\A_{i}$
and $\D$ $F^{\M}_{\l}$-primary over $\cup_{i\le n}\B_{i}$.
It is easy to see that $\C_{n}$ is $s$-primary over $\cup_{i\le n}\A_{i}$
and so we may choose $\D$ so that
$\C_{n}\subseteq\D$. Choose $a'\in\D$ so that
$t(a',\C_{n})=t(a,\C_{n})$. Let $\A'$ be $s$-primary over
$\A\cup a'$.

{\bf Claim 1.} $\A'\da_{\A}\B$.

\proof Immediate by Lemma 1.4. $\eop$ Claim 1.

{\bf Claim 2.} For all $i\le n$, $\A'\da_{\A}\B_{i}$.

\proof Clearly it is enough to show that
$a'\da_{\A}\B\cup\A_{i}$. Let $I=\{ j\le n\vert\ j\ne i\}$.
By Claim 1 and (vi) above,

(*) $a'\da_{\C_{n}}\B$.

\noindent
By the choice of $\B$,
$\cup_{j\in I}\A_{j}\da_{\A_{i}}\B$ and so $\C_{n}\da_{\A_{i}}\B$.
With (*) above, this implies $a'\da_{\A_{i}}\B$.
Since $a'\da_{\A}\A_{i}$, $a'\da_{\A}\B\cup\A_{i}$. $\eop$ Claim 2.

Since $Dp(a,\A )>0$, there is $b\not\in\A'$ such that
$t(b,\A')$ is orthogonal to $\A$ and $b\da_{\A'}\D$.
By Claim 2 and [HS1] Corollary 4.8, $t(b,\D)$ is orthogonal
to $\B_{i}$ for all $i\le n$. It is easy to see that
this contradicts the assumption that $\M$ does not have $\l$-dop.
$\eop$

\th 1.6 Lemma. Assume that $\M$ is superstable,
$\A\subseteq\B$ are $s$-saturated and $\A\ne\B$.
Then there is a singleton $a\in\B-\A$ such that $t(a,\A )$ is regular.

\proof As in the case of superstable theories. $\eop$

\chapter{2. Superstable with ndop}

Throughout this section we assume that $\M$ is superstable and
does not have $\l (\M )$-dop.

\th 2.1 Definition.

(i) We say that $(P,f,g)=((P,<),f,g)$
is an $s$-free tree of ($s$-saturated) model
$\A$
if the following holds:

(i) $(P,<)$ is a tree without branches of length $>\o$,
$f:(P-\{ r\} )\rightarrow\A$ and $g:P\rightarrow P(\A )$,
where $r\in P$ is the root of $P$ and $P(\A )$ is the power set
of $\A$,

(ii) $g(r)$ is $s$-primary model (over $\empty$ i.e. saturated
model of power $\l (\M )$),

(iii)  if $t$ is not the root and $u^{-}=t$ then
$t(f(u),g(t))$ is orthogonal to $g(t^{-})$,

(iv) if $t=u^{-}$ then $g(u)$ is $s$-primary over $g(t)\cup f(u)$,

(v) Assume $T,V\subseteq P$ and $u\in P$ are such that

(a) for
all $t\in T$, $t$ is comparable with $u$,

(b) $T$ is downwards closed.

(c) if $v\in V$ then for all $t$ such that $v\ge t>u$, $t\not\in T$.

\noindent
Then
$$\bigcup_{t\in T}g(t)\da_{g(u)}\bigcup_{v\in V}g(v).$$

(ii) We say that $(P,f,g)$ is
an $s$-decomposition of $\A$ if it is a maximal
$s$-free tree of
$\A$.

Notice that by Lemma 0.2 (iii) it is easy to see, that
every $e$-saturated model
has an $s$-decomposition.

\th 2.2 Theorem. ($\M$ superstable without $\l (\M )$-dop) Assume $\A$
is $e$-satu\-rated and $(P,f,g)$ is an $s$-decomposition
of $\A$. If $\B\subseteq\A$ is $s$-primary over
$\cup_{t\in P}g(t)$, then $\B=\A$.

\proof Immediate by Lemma 0.2 (iii) and (the proof of)
[HS1] Theorem 5.13. $\eop$

\th 2.3 Corollary. Suppose $\A$ and $\B$ are $e$-saturated.
If $(P,f,g)$ is a decomposition of
both $\A$ and $\B$, then $\A\cong\B$.

\proof Easy by Theorem 2.2. $\eop$

We say that an $s$-free tree $(P,f,g)$ is regular
if the following holds: if $t,u\in P$ are such that $u$ is an immediate
successor of $t$, then
$t(f(u),g(t))$ is regular.
We say that $(P,f,g)$ is a regular $s$-decomposition of $e$-saturated $\A$,
if it an $s$-decomposition of $\A$ and a regular $s$-free tree.

\th 2.4 Lemma. Every $e$-saturated model has a regular
$s$-decomposition.

\proof Immediate by Lemmas 0.2 (iii) and 1.6. $\eop$

\th 2.5 Definition.

(i) We say that $\M$ is shallow if every branch
in every regular $s$-free tree is finite.
If $\M$ is not shallow, then
we say that $\M$ is deep.

(ii) If $P=(P,<)$ is a tree without infinite branches, then
by $Dp(P)$ we mean the depth of $P$.

(iii) Assume that $\M$ is shallow.
We define the depth of $\M$ to be
$$sup\{ Dp(P)+1\vert\ (P,f,g)\ \hbox{\sl is a regular $s$-free tree}\} .$$

\th 2.6 Lemma. Assume that $\M$ is shallow and $\l (\M )$ is regular.
Then the depth of
$\M$ is $<\l (\M )^{+}$.

\proof Choose a minimal regular $s$-free tree $(P,f,g)$ so that if
$t\in P$ and $p\in S(g(t))$ are such that
if $t$ has an immediate predecessor $t^{-}$, then
$p$ is orthogonal to
$g(t^{-})$, then
there is an immediate
successor $u\in P$ of $t$ such that $t(f(u),g(t))=p$.
Clearly $Dp(P)<\l (\M )^{+}$. Also if $(P',f',g')$ is a
regular $s$-free tree, then there is an order-preserving function
$h:P'\rightarrow P$. Then $Dp(P')\le Dp(P)$, from which
the claim follows. $\eop$

By $\vert L\vert$ we mean the number of $L$-formulas modulo
the equivalence relation
$\models\forall x(\phi (x)\leftrightarrow\psi (x))$.

\th 2.7 Theorem. Assume that $\M$ is shallow. Then the depth of
$\M$ is $<(\vert S(\empty )\vert^{\o})^{+}$ and so
it is $<(2^{\vert L\vert})^{+}$.

\proof By Lemma 2.6, we may assume that $\l (\M )>\o$.
Choose a minimal regular $s$-free tree $(P,f,g)$ so that if
$t\in P$ and $p\in S(g(t))$ are such that
if $t$ has an immediate predecessor $t^{-}$, then
$p$ is orthogonal to
$g(t^{-})$, then
there is an immediate
successor $u\in P$ of $t$ and an automorphism $h$ of
$g(t)$ such that
such that $t(f(u),g(t))=h(p)$.

{\bf Claim 1.} $Dp(P)<(\vert S(\empty )\vert^{\o})^{+}$.

\proof Clearly it is enough to show that for all $t\in P$
the number of immediate successors of $t$ is
at most $\vert S(\empty )\vert^{\o}$. As in the proof of
Lemma 0.2, for all $p\in S(g(t))$, there is a countable
indiscernible $I\subseteq g(t)$ such that $Av(I,g(t))=p$.
Also if $t(I,\empty )=t(I',\empty )$, then there is
an automorphism $h$ of $g(t)$ such that
$h(I)=I'$ (remember that $g(t)$ is an
$F^{\M}_{\vert g(t)\vert}$-saturated model of
power $\l (\M )>\o$). So the number
of immediate successors of $t$ is at most
$$\vert\{ t(I,\empty )
\vert\ I\subseteq g(t)\ \hbox{\rm countable indiscernible}\}\vert .$$
Clearly this is at most
$\vert S(\empty )\vert^{\o}$. $\eop$ Claim 1.

{\bf Claim 2.} If $(P',f',g')$ is a
regular $s$-free tree, then there is an order-preserving function
$h:P'\rightarrow P$.

\proof Just choose $h$ so that

(i) if $r$ is the root of $P'$ then $h(r)$ is the root of $P$,

(ii) if $t'\in P'$ is not a root of $P'$ and
$u'$ is the immediate predecessor of $t'$, then $t=h(t')$
is such that it is an immediate successor of $u=h(u')$ and
there is an isomorphism $h^{*}:g'(u')\rightarrow g(u)$
satisfying $t(f(t),g(u))=h^{*}(t(f'(t'),g'(u')))$.

\noindent
Clearly this is possible. $\eop$ Claim 2.

As in the proof of Lemma 2.6, Claim 1 and 2 imply that
the depth of $\M$ is $<(\vert S(\empty )\vert^{\o})^{+}$. $\eop$

\th 2.8 Theorem. Assume that $\M$ is shallow and $\g^{*}$ is the
depth of $\M$. Then the number of
non-isomorphic $e$-saturated models
of power $\aleph_{\a}$ is at most
$\beth_{\g^{*}}(\vert\a\vert +\l (\M ))$.

\proof By Corollary 2.3, it is enough to count the number
of 'non-isomorphic' regular $s$-free trees
$(P,f,g)$ of power $\aleph_{\a}$.
This is an easy induction on $Dp(P)$, see the related results
in [Sh3]. $\eop$

\th 2.9 Theorem. Assume that $\M$ is shallow and $\g^{*}$ is
the depth of $\M$. Let $\k =\beth_{\g^{*}}(\l (\M ))^{+}$.
If $\A_{i}$, $i<\k$, are $e$-saturated models, then there are
$i<j<\k$ such that $\A_{i}$ is elementarily embeddable into $\A_{j}$.

\proof By Theorem 2.2, this question can be reduced to the
question of 'embeddality' of labelled trees. So this follows
immediately from [Sh3] X Theorem 5.16C. $\eop$

A cardinal $\k$ is called beautiful if $\k =\o$ or
for all $\xi <\k$,
$\k\warrow (\o )^{<\o}_{\xi}$, see [Sh2] Definition 2.3.

\th 2.10 Theorem. ($\M$ is superstable without $\l (\M )$-dop
but not necessarily
shallow.) Assume that there is a beautiful cardinal $>\l (\M )$.
Let $\k^{*}$ be the least such cardinal.
If $\A_{i}$, $i<\k^{*}$, are $e$-saturated models, then there are
$i<j<\k^{*}$ such that $\A_{i}$ is elementarily embeddable into $\A_{j}$.

\proof Again by Theorem 2.2, this follows immediately
from [Sh2] Theorems 5.8 and 2.10. $\eop$

If $(P,<)$ is a tree without branches of length $\ge\o$ and
$t\in P$, then by $Dp(t,P)$ we mean the depth of $t$ in $P$.
If $t$ is not the root, then
by $t^{-}$ we mean the immediate predecessor of $t$.

\th 2.11 Theorem. Assume that $\M$ is superstable,
deep, does not have $\l (\M )$-dop
and $(\l (\M ))^{+}$-dop and $\l>\l (\M )$. Then there are
$s$-saturated (and so $e$-saturated)
models $\A_{i}$, $i<2^{\l}$, of power $\l$ such that
for all $i<j<2^{\l}$, $\A_{i}\not\cong\A_{j}$.

{\bf Remark.} Assume $\M$ is superstable.
In the next section we show that $\M$ has many
$e$-saturated models
if $\M$ has $\l (\M )$-dop. Similarly
we can show that $\M$ has many
$e$-saturated models if $\M$ has
$(\l (\M ))^{+}$-dop. In fact,
it can be seen that $\l (\M )$-ndop
implies $(\l (\M ))^{+}$-ndop ($\l (\M )$-ndop
implies structure theorem for $s$-saturated and so especially for
$F^{\M}_{(\l (\M ))^{+}}$-saturated models, while
$(\l (\M ))^{+}$-dop implies a lot of non-structure for
$F^{\M}_{(\l (\M ))^{+}}$-saturated models).

\proof Assume $X_{i}\subseteq\l$, $i<2$, are such that
$X_{0}\ne X_{1}$. Choose
regular $s$-free trees $(P_{i},f_{i},g_{i})$, $i<2$, so that

(i) $P_{i}$ does not have branches of length $\ge\o$ but for all
$t\in P_{i}$, if $t$ is not the root, then $Dp(f(t),g(t^{-}))>0$,

(ii) for all $\a\in X_{i}$, there are $\l$ many $t\in P_{i}$
such that the height of $t$ is one and $Dp(t,P_{i})=\a$ and
if $Dp(t,P_{i})=\b$ and the height of $t$ is one, then $\b\in X_{i}$,

(iii) for all $t\in P_{i}$, if $Dp(t,P_{i})=\a$ and $\b <\a$, then
$\vert\{ u\in P_{i}
\vert\ u^{-}=t\ \hbox{\rm and}\ Dp(u,P_{i})\ge\b\}\vert =\l$,

(iv) if $t,u\in P_{i}$ are not the root and
$t^{-}=u^{-}$, then
$$t(f_{i}(t),g_{i}(t^{-}))=t(f_{i}(u),g_{i}(u^{-})),$$
we write $p_{t^{-}}$ for this type.

\noindent
Let $r_{i}$ be the root of $P_{i}$,
Choose finite $A_{i}\subseteq B_{i}\subseteq g_{i}(r_{i})$ so that
$p_{r_{i}}$ does not split strongly over $A_{i}$ and
$(p_{r_{i}}\raj B_{i},A_{i})$ is a regular stationary pair.
Then we require also

(v) $B_{0}=B_{1}$ (=$B$), $A_{0}=A_{1}$ (=$A$) and
$p_{r_{0}}\raj B=p_{r_{1}}\raj B$.

Let $\A_{i}$, $i<2$, be $s$-primary over $\cup_{t\in P_{i}}g_{i}(t)$.
We show that there is no isomorphism
$F:\A_{0}\rightarrow\A_{1}$ such that $F\raj B=id_{B}$.
Clearly this is enough (since $\l^{<\o}<2^{\l}$, 'naming'
finite number of elements does not change the number of models
and since $\M$ is $\l$-stable, $\vert\A_{i}\vert =\l$).
For a contradiction
we assume that $F$ exists. Clearly we may assume that $F=id_{\A_{0}}$,
this simplifies the notation.

We let $P^{*}_{i}$ be the set of those $t\in P_{i}$, which are
not leafs.
For all $t\in P^{*}_{0}$, we let $G(t)\in P^{*}_{1}$ be
(some node) such that
$p_{t}$ is not orthogonal to $p_{G(t)}$ (if exists).

{\bf Claim.} $G$ is an one-to-one function from $P^{*}_{0}$ onto
$P^{*}_{1}$.

\proof Since for all $t\in P^{*}_{0}$,
$\vert\{ u\in P_{0}\vert\ u^{-}=t\}\vert =\l >\l (\M )$,
the existence of $G(t)$ follows easily.
Since for all $u,u'\in P^{*}_{1}$, $u\ne u'$, $p_{u}$ is orthogonal to
$p_{u'}$, $G(t)$ is unique by Corollary 1.3. But then by symmetry,
claim follows. $\eop$ Claim.

We prove a contradiction (with (i) above) by constructing a
strictly increasing sequence $(t_{j})_{j<\o}$ of elements of
$P^{*}_{0}$. We construct also a strictly increasing sequence
$(u_{j})_{j<\o}$ of elements of $P_{1}$, sets $I^{i}_{j}$, $i<2$,
and models $\B_{j}$ so that

(1) $Dp(u_{j},P_{1})<Dp(t_{j},P_{0})$ and for all
$t\ge t_{j}$, $G(t)\ge u_{j}$,

(2) $I^{i}_{j}\subseteq P_{i}$ is downwards closed,
non-empty and of power
$\le\l (\M )$ and $I^{i}_{j}\subseteq I^{i}_{j+1}$,

(3) $t_{j}\in I^{0}_{j+1}$ and $G(t_{j})\in I^{1}_{j+1}$,

(4) $\B_{j}$ is $s$-primary over $\cup_{t\in I^{0}_{j}}g_{0}(t)$
and over $\cup_{u\in I^{1}_{j}}g_{1}(u)$ and $\B_{j}\subseteq\B_{j+1}$.

\noindent
We do this by induction on $j<\o$.

$j=0$: Choose $I^{0}_{0}$, $I^{1}_{0}$ and $\B_{0}$ so that
(2) and (4) above are satisfied (if $\B'\subseteq\B_{0}$
is $s$-primary over $\cup_{t\in I}g(t)$, $I\subseteq P_{0}$, then by Theorem
2.2 and [HS1] Lemma 5.4 (ii), $\B_{0}$ is $s$-primary
over $\B'\cup\bigcup_{t\in P_{0}}g(t)$).
Let $t_{0}\in P_{0}$ be such that $t_{0}\not\in I^{0}_{0}$
and $(t_{0})^{-}=r_{0}$. Then
$$(*)\ \ \ \ f_{0}(t_{0})\da_{A}\B_{0}.$$
By Lemma 1.5, there is $u_{0}\in P_{1}-I^{1}_{1}$ such that
$f_{1}(u_{0})\nda_{\B_{0}}f_{0}(t_{0})$ and
$(u_{0})^{-}\in I^{1}_{0}$. By Lemma 1.4,
$$f_{0}(t_{0})\da_{\B_{0}}\cup\{g_{1}(u)\vert\ u\not\ge u_{0}\} .$$
So $u_{0}$ is unique and
the latter half of (1) holds. By (*), $(u_{0})^{-}=r_{1}$
and so since $X_{0}\ne X_{1}$ we can choose $t_{0}$ so that
$Dp(u_{0},P_{1})\ne Dp(t_{0},P_{0})$. By symmetry, we may assume that
$Dp(u_{0},P_{1})<Dp(t_{0},P_{0})$. Finally, this implies that
$t_{0}\in P^{*}_{0}$.

$j=k+1$: Essentially, just repeat the argument above.
$\eop$

\chapter{3. Superstable with dop or unstable}

\th 3.1 Theorem. Assume $\M$ is superstable with $\l (\M )$-dop,
$\k >(\l_{r}(\M ))^{+}$ is regular and $\xi >\k$.
Then there are $F^{\M}_{\k}$-saturated (and so $e$-saturated)
models $\A_{i}$, $i<2^{\xi}$, of power $\xi$ such that for all
$i\ne j$, $\A_{i}$ is not elementarily embeddable into $\A_{j}$.

\proof By [HS1] Corollary 6.5 and (the proof of)
[Hy] Lemma 2.5,
this follows from [Sh4] Theorems 3.20 and 3.27
and the claim below:
Let a linear ordering $\n$ be almost $\k$-homogeneous i.e.
for all $X\subseteq\n$ of power $<\k$ there is $Y\subseteq\n$
of power $<\k$ such that $X\subseteq Y$ and if $x,y\in\n$ are
in the same Dedekind cut of $Y$, then there is an automorphism
$f$ of $\n$ such that $f\raj Y=id_{Y}$ and $f(x)=y$.
Let $\A_{\n}$, $\phi$, $\psi$ and $B_{i}$, $C_{i}$ and $I_{ij}$, $i,j\in\n$
be as in [Hy]. For all $X\subseteq\n$, by
$S_{X}$ we mean the set
$\bigcup\{ B_{i}\cup C_{i}\vert\ i\in X\}\cup\bigcup
\{ I_{i,j}\vert\ i,j\in X ,\ i<j\}$.

{\bf Claim.} $(B_{i}\cup C_{i})_{i\in\n}$ is weakly
$(\k ,\phi )$-skeleton-like in $\A_{\n}$ (see [Sh4]).

\proof Let $A\subseteq\A_{\n}$ be of power $(\l_{r}(\M ))^{+}$.
Since $\k$ is regular, we can find $X\subseteq\n$ of power $<\k$
and $B\subseteq\A_{\n}$ of power $<\k$
such that

(i) $A\subseteq B$,

(ii) $\A_{\n}$ is $F^{\M}_{\k}$-primary over $B\cup S_{\n}$,

(iii) for all $a\in B$,
$t(a,S_{\n})\in F^{\M}_{\k}(S_{X})$.

\noindent
Let $Y\subseteq\n$ be as in the definition of almost $\k$-homogeneous.
Let $x,y\in\n$ be in the same Dedekind cut of $Y$ and
assume that $\A_{\n}\models\psi (A,B_{x}\cup C_{y})$.
It is
enough to show that $\A_{\n}\models\psi (A,B_{y}\cup C_{y})$.

By the choice of $Y$, there is an automorphism
$f$ of $\n$ such that $f\raj Y=id_{Y}$ and $f(x)=y$.
This $f$ induces an elementary function $g$ from $S_{\n}$ onto $S_{\n}$
such that $g\raj S_{Y}=id_{S_{Y}}$ and $g\raj B_{x}\cup C_{x}$
is the natural elementary function onto $B_{y}\cup C_{y}$.
By (iii) above, we can find an automorphism $h$ of $\M$
such that $g\subseteq h$ and $h\raj B=id_{B}$.
Let $\A'_{\n}=h(\A_{\n})$. Then both of the models are
$F^{\M}_{\k}$-primary over $B\cup S_{\n}$ and so they
are isomorphic over $B\cup S_{\n}$. Let
$h'$ be the isomorphism from $\A'_{\n}$ to $\A_{\n}$.
Then $h'\circ h\raj\A_{\n}$ is
an automorphism of $\A_{\n}$, $h'\circ h\raj A=id_{A}$
and $h'\circ h\raj B_{x}\cup C_{x}$
is the natural elementary function onto $B_{y}\cup C_{y}$.
Clearly this implies that $\A_{\n}\models\psi (A,B_{y}\cup C_{y})$.
$\eop$ Claim.

$\eop$

\th 3.2 Lemma. Assume that $\M$ is unstable.
Let $\k >\vert L\vert$ be
a regular cardinal,
and $\n=(\n ,<)$ be a linear ordering. Then there are
sequences $a_{i}$, $i\in\n$, a model $\A$ and functions
$f_{i}:\M^{n_{i}}\rightarrow\M$, $i<2^{<\k}$,
such that $n_{i}<\o$ and if we write
$L^{*}=L\cup\{ f_{i}\vert\ i<2^{<\k}\}$
then the following holds:

(i) $(a_{i})_{i\in\n}$ is order-indiscernible
inside $\A$ in the language $L^{*}$,

(ii) for all $X\subseteq\n$, the closure $\A_{X}$ of $\{ a_{i}\vert\ i\in X\}$
under the functions of $L^{*}$ is a locally $F^{\M}_{\k}$-saturated model
(in the language $L$) and $\A=\A_{\n}$,

(iii) there is an $L$-formula $\phi (x,y)$
such that for all $i,j\in\n$, $\models\phi (a_{i},a_{j})$ iff $i<j$.

\proof Define functions $f'_{i}:\M^{n_{i}}\rightarrow\M$, $i<2^{<\k}$,
so that

(*) the closure of any set
under the functions $f_{i}$ is locally
$F^{\M}_{\k}$-saturated (in $L$) and
$L'$-elementary submodel of
$(\M ,f'_{i})_{i<2^{<\k}}$, where $L'=L\cup\{ f'_{i}\vert\ i<2^{<\k}\}$.

By Erd\"os-Rado Theorem and [Sh1] I Lemma 2.10 (1),
we can find sequences $(a^{k}_{i})_{i<k}$, $k<\o$, such that

(1) there is a formula $\phi (x,y)$ such that for all
$k<\o$ and $i,j<k$,
$\models\phi (a^{k}_{i},a^{k}_{j})$ iff $i<j$,

(2) $(a^{k}_{i})_{i<k}$ is order-indiscernible in the language $L'$,

(3) the $L'$-type of $(a^{k}_{i})_{i<k}$ (over $\empty$)
is the same as the $L'$-type of $(a^{k+1}_{i})_{i<k}$.

\noindent
Since $\M$ is homogeneous, we can find for all $i\in\n$, $a_{i}$
so that for all $k<\o$,
if $i_{0}<i_{1}<...<i_{k-1}$, then
$t((a_{i_{j}})_{j<k},\empty )=t((a^{k}_{j})_{j<k},\empty )$.
Again, since $\M$ is homogeneous (use e.g. [HS1] Lemma 1.1)
we can define the functions $f_{i}$ so that
for all $i_{0}<i_{1}<...<i_{k-1}$
the following holds:

(**) If $\A_{1}$ is the closure of
$(a_{i_{j}})_{j<k}$ under the functions $f_{i}$
and $\A_{2}$ is
the closure of $(a^{k}_{j})_{j<k}$
under the functions $f'_{i}$,
then there is an $L$-isomorphism $F:\A_{1}\rightarrow\A_{2}$,
such that $F(a_{i_{j}})=a^{k}_{j}$
and for all $a,b\in\A_{1}$ and $i<2^{<\k}$, $f_{i}(a)=b$
iff $f'_{i}(F(a))=F(b)$.

\noindent
Let $\A=\A_{\n}$, i.e. the closure of $\{ a_{i}\vert\ i\in\n\}$
under the functions of $L^{*}$.
Then it is easy to see that (iii) in
the claim is satisfied.

(ii): Assume $X\subseteq\n$. We show
that $\A_{X}$ is locally $F^{\M}_{\k}$-saturated. For this
let $A\subseteq\A_{X}$ be finite. Then there is $X'\subseteq X$
finite, such that $A\subseteq\A_{X'}$. By (**) above,
$\A_{X'}$ is locally $F^{\M}_{\k}$-saturated. So there is
$F^{\M}_{\k}$-saturated $\B$ such that $A\subseteq\B\subseteq\A_{X}$.

(i): By (*) and (**) above
it is easy to see that for all finite $X\subseteq\n$,
$\A_{X}$ is an $L^{*}$-elementary submodel of $\A$. By (2), (*) and (**)
again, (i) follows.
$\eop$

\th 3.3 Theorem. Assume $\M$ is unstable.
Let $\l$ and $\k$ be
regular cardinals, $\l >2^{<\k}$ and $\k >\vert L\vert$.
Then there are locally $F^{\M}_{\k}$-saturated models
$\A_{i}$, $i<2^{\l}$, such that $\vert\A_{i}\vert =\l$ and
if $i\ne j$, then
$\A_{i}$ is not elementarily embeddable into $\A_{j}$.

\proof By Lemma 3.2 this follows from [Sh4] Chapter 6 Theorem 3.1 (3).
Notice that the trees can be coded into linear orderings. $\eop$

\chapter{4. Strictly stable}

Through out this section we assume that $\M$ is stable but unsuperstable,
and that $\k =cf(\k )>\l_{r}(\M )$.

We write $\k^{\le\o}$ for
$\{ \n :\a\rightarrow\k\vert\ \a\le\o\}$,
$\k^{<\o}$ and $\k^{\o}=\k^{=\o}$ are defined similarly
(of course these have also the other meaning, but it will
be clear from the context, which one we mean).
Let $J\subseteq 2^{\le\k}$.
We order $P_{\o}(J)$ (=the set
of all finite subsets of $J$) by defining
$u\le v$ if for every $\n\in u$ there is $\xi\in v$ such that
$\n$ is an initial segment of $\xi$.

Since $\M$ is unsuperstable, by [HS1] Lemma 5.1, there are
$a$ and $F^{\M}_{\l_{r}(\M )}$-satu\-rated models $\A_{i}$, $i<\o$,
of power $\l_{r}(\M )$ such that

(i) if $j<i<\o$, then $\A_{j}\subseteq\A_{i}$,

(ii) for all $i<\o$, $a\nda_{\A_{i}}\A_{i+1}$.

\noindent
Let $\A_{\o}$ be an $F^{\M}_{\l_{r}(\M )}$-primary model
over $a\cup\bigcup_{i<\o}\A_{i}$.
Then for all $\n\in\k^{\le\o}$,
we can find $\A_{\n}$ such that

(a) for all $\n\in\k^{\le\o}$, there is an automorphism $f_{\n}$ of $\M$
such that $f_{\n}(\A_{length(\n )})$ $=\A_{\n}$,

(b) if $\n$ is an initial segment of $\xi$, then
$f_{\xi}\raj\A_{length(\n )}=f_{\n}\raj\A_{length(\n )}$,

(c) if $\n\in\k^{<\o}$, $\a\in\k$ and
$X$ is the set of those $\xi\in\k^{\le\o}$ such that
$\n\frown(\a )$ is an initial segment of $\xi$, then
$$\cup_{\xi\in X}\A_{\xi}\da_{\A_{\n}}\cup_{\xi\in(\k^{\le\o}-X)}\A_{\xi}.$$

\noindent
For all $\n\in\k^{\o}$, we let $a_{\n}=f_{\n}(a)$.

For each $\a <\k$ of cofinality $\o$, let $\n_{\a}\in\k^{\o}$
be a strictly increasing sequence such that $\cup_{i<\o}\n_{\a}(i)=\a$.
Let $S\subseteq\{\a <\k\vert\ cf(\a )=\o\}$.
By $J_{S}$ we mean the set
$$\k^{<\o}\cup\{ \n_{\a}
\vert\ \a\in S\} .$$
Let $I_{S}=P_{\o}(J_{S})$.

\th 4.1 Lemma. For all $S\subseteq\{\a <\k\vert\ cf(\a )=\o\}$,
there are sets $\A_{u}$, $u\in I_{S}$, such that

(i) for all $u,v\in I_{S}$, $u\le v$ implies $\A_{u}\subseteq\A_{v}$,

(ii) for all $u\in I_{S}$, $\A_{u}$ is $F^{\M}_{\l_{r}(\M )}$-primary over
$\cup_{\n\in u}A_{\n}$,

(iii) if $\a\in\k -S$, $u\in I_{S}$
and $v\in P_{\o}(J_{S}\cap\a^{\le\o})$ is maximal
such that $v\le u$, then
$$\A_{u}\da_{\A_{v}}\cup_{w\in P_{\o}(J_{S}\cap\a^{\le\o})}\A_{w}.$$

\proof See [HS2] Lemmas 4 and 7. $\eop$

For all $S\subseteq\{\a <\k\vert\ cf(\a )=\o\}$, let
$\A_{S}=\cup_{u\in I_{S}}\A_{u}$. By Lemma 4.1 (i) and (ii),
$\A_{S}$ is $e$-saturated and $\vert\A_{S}\vert =\k$.

\th 4.2 Lemma. There are sets
$S_{i}\subseteq\{\a <\k\vert\ cf(\a )=\o\}$, $i<2^{\k}$,
such that if $i\ne j$, then $S_{i}-S_{j}$ is stationary.

\proof Let $f_{i};\k\rightarrow\k$, $i<2$, be one to one functions
such that $rng(f_{0})\cap rng(f_{1})=\empty$. Let
$R'_{i}$, $i<2^{\k}$, be an enumeration of the power set
of $\k$. We define $R_{i}$, $i<2^{\k}$, so that
$f_{0}(\a )\in R_{i}$ iff $\a\in R'_{i}$ and
$f_{1}(\a )\in R_{i}$ iff $\a\not\in R'_{i}$. Then clearly,
$i\ne j$ implies $R_{i}-R_{j}\ne\empty$. By [Sh3] Appendix Theorem 1.3 (2),
there are pairwise disjoint
stationary sets $S'_{j}\subseteq\{\a <\k\vert\ cf(\a )=\o\}$, $j<\k$.
For $i<2^{\k}$,
we let $S_{i}=\cup_{j\in R_{i}}S'_{j}$. Clearly these are
as wanted. $\eop$

\th 4.3 Theorem. Assume $\M$ is stable and unsuperstable and
$\k =cf(\k )>\l_{r}(\M )$. Then there are $e$-saturated
models $\A_{i}$, $i<2^{\k}$, of power $\k$
such that if $i\ne j$,
then $\A_{i}$ is not elementarily embeddable into $\A_{j}$.

\proof For all $i<2^{\k}$, let $\A_{i}=\A_{S_{i}}$, where
the sets $S_{i}$
are as in Lemma 4.2. Assume $i\ne j$. We show that
there are no elementary map $F:\A_{i}\rightarrow\A_{j}$.

For a contradiction, assume that $F$ exists. For all $\a<\k$, let
$I^{\a}_{S_{i}}$ be the set of those $u\in I_{S_{i}}$ such that
for all $\n\in u$, $sup\{\n (i)\vert\ i<length(\n )\} <\a$.
Let $\A^{\a}_{i}=\cup_{u\in I^{\a}_{S_{i}}}\A_{u}$.
$I^{\a}_{S_{j}}$ and $\A^{\a}_{j}$ are defined similarly.
We say that $\a$ is closed if for all $a\in\A_{i}$,
$a\in\A^{\a}_{i}$ iff $F(a)\in\A^{\a}_{j}$.
Let $C$ be the set of all closed ordinals and
$C_{lim}$ the set of all limit points in $C$.
Then $S^{0}=C_{lim}\cap (S_{i}-S_{j})$ is stationary.

For all $\a\in S^{0}$, let $u_{\a}\in I_{S_{j}}$ be such that
$F(a_{\n_{\a}})\in\A_{u_{\a}}$. By $g(\a )$ we mean the
least $\b\in C$ such that
$u_{\a}\da_{\A^{\b}_{j}}\A^{\a}_{j}$.
By Lemma 4.1 (iii) and the fact that
$S^{0}\cap S_{j}=\empty$, $g(\a )<\a$. So there is stationary
$S^{1}\subseteq S^{0}$ such that $g\raj S^{1}$ is constant.
Let $\a^{*}$ be this constant value.

Then there is $S^{2}\subseteq S^{1}$ and $n<\o$ such that
$\vert S^{2}\vert =\k$ and for all $\b ,\g\in S^{2}$,
if $\b\ne\g$, then $\n_{\b}(n)\ne\n_{\g}(n)$. By choosing $n$ so that
it is minimal, we may assume that for all $\b\in S^{2}$,
$\n_{\b}(n-1)<\a^{*}$. Clearly we may assume that for all
$\b\in S^{2}$, $\n_{\b}(n)>\a^{*}$.

Then by Lemma 4.1 (iii),

(i) $(F(\A_{\n_{\b}\raj (n+1)}))_{\b\in S^{2}}$ is
$F(\A^{\a^{*}}_{i})$-independent.

\noindent
Since $F(a_{\n_{\b}})\da_{\A^{\a^{*}}_{j}}F(\A_{\n_{\b}\raj (n+1)})$
and $F(a_{\n_{\b}})\nda_{F(\A^{\a^{*}}_{i})}F(\A_{\n_{\b}\raj (n+1)})$,

(ii) for all $\b\in S^{2}$,
$F(\A_{\n_{\b}\raj (n+1)})\nda_{F(\A^{\a^{*}}_{i})}\A^{\a^{*}}_{j}$.

\noindent
Since $\k (\M )<\k$, $\vert \A^{\a^{*}}_{j}\vert <\k$
and $\vert S^{2}\vert =\k$, (i) and (ii) are contradictory.
$\eop$

\vfill
\eject

\chapter{References}

\item{[Hy]} T. Hyttinen, On nonstructure of elementary submodels
of a stable homogeneous structure, Fundamenta Mathematicae, to appear.

\item{[HS1]} T. Hyttinen and S. Shelah, Strong splitting in stable
homogeneous models, submitted.

\item{[HS2]} T. Hyttinen and S. Shelah, On the number of elementary
submodels of an unsuperstable homogeneous structure, Mathematical
Logic Quarterly, to appear.

\item{[Sh1]} S. Shelah, Finite diagrams stable in power, Annals
of Mathematical Logic, vol. 2, 1970, 69-118.

\item{[Sh2]} S. Shelah, Better Quasi-orders for uncountable
cardinals, Israel Journal of Mathematics, vol. 42, 1982, 177-226.

\item{[Sh3]} S. Shelah, Classification Theory, Stud. Logic Found. Math.
92, North-Holland, Amsterdam, 2nd rev. ed., 1990.

\item{[Sh4]} S. Shelah, Non-structure Theory, to appear.

\bigskip

\settabs\+\hskip 0truecm&\hskip 6truecm&\cr

\+Tapani Hyttinen\cr

\+Department of Mathematics\cr

\+P.O. Box 4\cr

\+00014 University of Helsinki\cr

\+Finland\cr

\medskip

\+&Saharon Shelah\cr

\+&Institute of Mathematics&Rutgers University\cr

\+&The Hebrew University&Hill Ctr-Bush\cr

\+&Jerusalem&New Brunswick\cr

\+&Israel&New Jersey 08903\cr
\+&&U.S.A.\cr

\end